\documentclass{article}
\usepackage{amsmath,amssymb}
\usepackage{theorem}
\usepackage{color}

\theorembodyfont{\itshape}
\theoremstyle{plain}
\newtheorem{theorem}{Theorem}[section]
\newtheorem{lemma}{Lemma}[section]
\newtheorem{corollary}{Corollary}[section]
\newtheorem{proposition}{Proposition}[section]
\theorembodyfont{\rmfamily}
\newtheorem{example}{Example}[section]
\newtheorem{question}{Question}[section]
\newtheorem{problem}{Problem}[section]
\newtheorem{definition}{Definition}[section]

\newtheorem{remark}{Remark}[section]
\newtheorem{claim}{Claim}[section]

\newcommand{\Int}{\mbox{{\rm Int}}\,}
\newcommand{\Cl}{\mbox{{\rm Cl}}}

\begin{document}

\title{\bf On $\pi$-compatible topologies and their special cases}

\author{Vitalij A.~Chatyrko}


\maketitle

\begin{abstract} Topologies 
$\tau, \sigma$ on a set $X$ are called {\it $\pi$-compatible} if 
$\tau$ is a $\pi$-network for $\sigma$,  and vice versa. 
If topologies $\tau, \sigma$ on a set $X$ are $\pi$-compatible then the families of 
nowhere dense sets (resp. meager sets or sets possessing the Baire property) of the spaces $(X, \tau)$ and $(X, \sigma)$ coincide.
A topology $\sigma$ on a set $X$ is called {\it an admissible extension} of a topology $\tau$ on $X$
if $\tau \subseteq \sigma$ and $\tau$ is a $\pi$-network for $\sigma$. It turns out that 
examples of admissible extensions
were occurred in literature several times.
In the paper we provide some new facts about the $\pi$-compatibility and the admissible extension
as well as about their particular cases.
\end{abstract} 


\smallskip
{\it 2000 AMS (MOS) Subj. Class.:} Primary 54A10, 54D10, 54D15   
\medskip
\baselineskip=18pt

\section{Introduction}
Let $X$ be a non-empty set, $\mathcal{P}(X)$ the family of all subsets of $X$,  $\tau$ a topology on $X$ and  $\mathcal{I}_m(X, \tau)$ the family of all meager sets of the topological space $(X,\tau)$. 

An interesting collection of subsets of $X$ extending $\tau$ as well as $\mathcal{I}_m(X, \tau)$,  is the family $\mathcal{B}_p(X, \tau)$ of all subsets of $X$ possessing the Baire property in $(X,\tau)$. Recall that a subset $A$ of $X$ has the Baire property in the space $(X, \tau)$ if $A = (O \setminus M) \cup N$, where $O \in \tau$ and $M, N \in \mathcal{I}_m(X, \tau)$.  It is well known (cf. \cite{Ku}) that the family $\mathcal{B}_p(X,\tau)$ is a $\sigma$-algebra of sets which is  invariant under automorphisms of the space 
 $(X,\tau)$. 
  
Let us note that $\mathcal{B}_p(X,\tau_{\mbox{triv}}) = \{\emptyset, X\}$, $\mathcal{B}_p(X,\tau_{\mbox{dis}}) = \mathcal{P}(X)$, where 
$\tau_{\mbox{triv}}$ (resp. $\tau_{\mbox{dis}}$) is the trivial (resp. discrete) topology on the set $X$. It is easy to see that  for any set $X$ with $|X| >1$ there is no topology $\tau$ distinct from $\tau_{\mbox{triv}}$ such that 
$\mathcal{B}_p(X, \tau) = \{\emptyset, X\}$. But there are simple examples of
 a set $X$ and a topology $\tau$ on $X$ distinct from 
$\tau_{\mbox{dis}}$ such that $\mathcal{B}_p(X,\tau) = \mathcal{P}(X)$ (see Corollary~\ref{remark1_baire}),  
and  a set $Y$ and a topology $\sigma$ on $Y$ distinct from $\tau_{\mbox{triv}}$ for which  $\mathcal{B}_p(Y,\sigma) \ne \mathcal{P}(Y)$ (see Example~\ref{topological sum} and Example~\ref{remark2_baire}).

Let us recall (cf. \cite{Ku}) that for the real numbers $\mathbb R$ with the Euclidean topology $\tau_E$ we have 
$\mathcal{B}_p(\mathbb R,\tau_E) \ne \mathcal{P}(\mathbb R)$, and  there is a lot of information 
about elements of the family $\mathcal{P}(\mathbb R) \setminus \mathcal{B}_p(\mathbb R,\tau_E)$ (see for example  \cite{Kh}). It would be interesting to know for what topologies $\tau$ on $\mathbb R$ the equality  $\mathcal{B}_p(\mathbb R,\tau ) = \mathcal{B}_p(\mathbb R,\tau_E)$ is valid.

One can pose  more general  questions.

\begin{question} \label{question 1}    Let $X$ be a set and $\tau$ be a topology on $X$.
\begin{itemize}
\item[(a)] Describe elements of the family $\mathcal{P}(X) \setminus \mathcal{B}_p(X,\tau)$.

\item[(b)] For what topologies  $\sigma$ on  the set $X$ does the equality 
$\mathcal{B}_p(X,\tau) = \mathcal{B}_p(X,\sigma)$ hold?
\end{itemize}
\end{question}

The $\pi$-compatibility  of topologies $\sigma$ and $\tau$ on a set $X$ 
(see Definition~\ref{pi compact}) introduced in \cite{CN1}
implies $\mathcal{I}_m(X,\tau) = \mathcal{I}_m(X,\sigma)$ (see Lemma~\ref{compatibility_meager}) and $\mathcal{B}_p(X,\tau) = \mathcal{B}_p(X,\sigma)$ (see 
Corollary~\ref{compatibility_baire}) but it is not equivalent to the equalities (see Example~\ref{remark_baire}).
It turns out that a stronger version of  $\pi$-compatibility between two topologies, the notion of the admissible
extension (see Definition~\ref{adm_ext}), was occurred in literature several times ( see \cite{H}, \cite{F}, \cite{CN2} and others). 
In the paper we provide some new facts about these relations between topologies which will be valid for the mentioned cases.

For standard notions we refer to \cite{Ku} and \cite{E}.

\section{Some simple answers to Question~\ref{question 1}.(a) }

\begin{proposition}\label{countable_pro} Let $X$ be a countable set, $Y \subseteq X$ and $\tau$ a topology on $X$.  Let also for each point $x \in Y$  the set $\{x\}$ is a  nowhere dense set in the space $(X, \tau)$, and 
for each point $x \in X \setminus Y$ the set $\{x\}$ is an open set of the space $(X, \tau)$. Then 
$\mathcal I_m(X,\tau) = \mathcal P(Y)$,  
$\mathcal B_p(X,\tau) = \mathcal P(X)$ and so $\mathcal P(X) \setminus \mathcal B_p(X,\tau) = \emptyset.$
$\Box$
\end{proposition}

\begin{corollary}\label{remark1_baire} 
\begin{itemize}
\item[(a)] Let $\mathbb Q$ be the set of rational numbers with the Euclidean topology $\tau_E$.
Then we have 
$\mathcal I_m(\mathbb Q,\tau_E) =  \mathcal B_p(\mathbb Q,\tau_E) = \mathcal P(\mathbb Q)$.

\item[(b)] Let $\mathbb Z$ be the set of integer numbers, $\mathbb Z_e$ be the set of even integer numbers and $\tau_K$ be a topology generated by the family $\{\{n-1, n, n+1 \}: n \in \mathbb Z_e\}$ (Khalimsky topology).
Then we have  $\mathcal I_m(\mathbb Z,\tau_K) = \mathcal P(\mathbb Z_e) \ne \mathcal P(\mathbb Z)= \mathcal B_p(\mathbb Z,\tau_K)$. $\Box$
\end{itemize}
\end{corollary}

The condition "each one-point set is either nowhere dense or open" in Proposition~\ref{countable_pro} is essential. In fact,

\begin{example}\label{topological sum} Let $\mathbb N$ be the set of positive integers and $\tau_p$ be the odd-even topology generated by the partition $\mathcal P = \{\{2k-1, 2k \} \}$ of $\mathbb N$  (\cite[Example 6]{SS}). Note that the family of nowhere dense sets (resp. meager sets) for the space  $(\mathbb N, \tau_p)$ is equal to $\{\emptyset\}$.
Moreover, $\mathcal{B}_p(\mathbb N, \tau_p) =  \tau_p \ne  \mathcal P(\mathbb N)$, and the family
$\mathcal P(\mathbb N) \setminus \mathcal{B}_p(\mathbb N, \tau_p)$ consists of sets which are not open.
\end{example}

Let us also note that the condition "to be countable" on the set $X$ in Proposition~\ref{countable_pro} is essential.  

\begin{example}\label{remark2_baire} 
Let $X$ be any uncountable set and $\tau_{\mbox{fc}}$ be the finite complement topology on $X$, i.e. $V \in \tau_{\mbox{fc}}$ iff $X \setminus V$ is finite. Since each infinite subset $Y$ of $X$ is dense in the space 
$(X, \tau_{\mbox{fc}})$, $\mathcal I_m(X, \tau_{\mbox{fc}}) = \{A \subset X: |A| \leq \aleph_0\}$ and 
$\mathcal B_p(X, \tau_{\mbox{fc}}) = \{A \subset X : |A| \leq \aleph_0 \mbox{ or } |X \setminus A| \leq \aleph_0\}.$ So any uncountable subset $Y$ of $X$ such that 
$X \setminus Y$ is also uncountable, does not belong to the family   $\mathcal{B}_p(X,\tau_{\mbox{fc}})$. In particular,  $\mathcal{B}_p(X,\tau_{\mbox{fc}}) \ne \mathcal P(X)$.
\end{example}

\section{$\pi$-compatible topologies as an answer to Question~\ref{question 1}.(b).}

\begin{definition}\label{pi compact} (\cite[Definition 3.1]{CN1}) Let $\tau, \sigma$ be topologies on a set $X$. We call the topologies {\it $\pi$-compatible} if 
$\tau$ is a $\pi$-network for $\sigma$ (i.e. for each non-empty element $O$ of $\sigma$ there is a non-empty element $V$ of $\tau$ which is a subset of $O$) and vice versa. 
\end{definition}

It is evident that the $\pi$-compatibility is an equivalence relation on the family of all topologies on the set $X$.

\begin{lemma}\label{compatibility_meager} (\cite[Lemma 3.3]{CN1}) Let $\tau$ and $\sigma$ be $\pi$-compatible topologies on a set $X$. Then the spaces $(X, \tau)$ and $(X, \sigma)$ have the same families of nowhere dense sets (respectively, meager sets). $\Box$
\end{lemma}

\begin{theorem}\label{compatibility_theroem} Let $\tau$ and $\sigma$ be $\pi$-compatible topologies on a set $X$. Then each non-empty element of $\tau$ is the union of 
a non-empty element of $\sigma$ and a nowhere dense set, and vice versa.
\end{theorem}

Proof. 
In fact, let $O_\tau$ be a non-empty element of $\tau$. Put $O_\sigma = \Int_\sigma  O_\tau$, where $\Int_\sigma A$ (respectively, $\Cl_\sigma A$) is  the interior (respectively, closure) of $A$ in the space $(X, \sigma)$ for any $A \subseteq X$. Note that $O_\sigma \ne \emptyset$ and $O_\tau = O_\sigma \cup ((\Cl_\sigma O_\sigma \setminus O_\sigma) \cap O_\tau) \cup (O_\tau \setminus \Cl_\sigma O_\sigma)$. Note also that the set $(\Cl_\sigma O_\sigma \setminus O_\sigma) \cap O_\tau$ is nowhere dense.

\begin{claim} The set $A = O_\tau \setminus \Cl_\sigma O_\sigma$ is nowhere dense.
\end{claim}
Proof.  Assume that the set $A$ is not nowhere dense. Put $W_\sigma = \Int_\sigma \Cl_\sigma  A$. Note that $W_\sigma \ne \emptyset$, $W_\sigma \subseteq \Cl_\sigma A$ and $W_\sigma \cap O_\sigma = \emptyset$. Since $W_\sigma \in \sigma$ and $\tau$ is a $\pi$-network for $\sigma$, there exists an non-empty element $T_\tau$ of $\tau$ such that $T_\tau \subseteq W_\sigma$.
There are two possibilities.

Case 1. $T_\tau \cap O_\tau \ne \emptyset$. Since $T_\tau \cap O_\tau \in \tau$ and $\sigma$ is a $\pi$-network for $\tau$,
there is a non-empty element $P_\sigma$ of $\sigma$ such that $P_\sigma \subseteq T_\tau \cap O_\tau$. Note that 
$P_\sigma \subseteq O_\tau \setminus O_\sigma$.
We have a contradiction with the definition of $O_\sigma$.

Case 2. $T_\tau \cap O_\tau = \emptyset$. Since $\sigma$ is a $\pi$-network of $\tau$, there exists a non-empty element $Q_\sigma$ of $\sigma$ such that $Q_\sigma \subseteq T_\tau$. So $Q_\sigma \cap O_\tau = \emptyset$ and hence 
$Q_\sigma \cap \Cl_\sigma O_\tau = \emptyset$. But $Q_\sigma \subseteq W_\sigma \subseteq \Cl_\sigma A \subseteq \Cl_\sigma O_\tau$. We have a contradiction.

So $A$ is nowhere dense. $\Box$

By the use of the claim we get that $O_\tau$ is  the union of $O_\sigma$ and the nowhere dense set  $((\Cl_\sigma O_\sigma \setminus O_\sigma) \cap O_\tau) \cup (O_\tau \setminus \Cl_\sigma O_\sigma)$. We have done. $\Box$

\begin{corollary}\label{compatibility_baire} (\cite[Theorem 3.4]{CN1}) Let $\tau$ and $\sigma$ be $\pi$-compatible topologies on a set $X$. Then the spaces $(X, \tau)$ and $(X, \sigma)$ have the same families of Baire sets. 
\end{corollary}
Proof. Let us consider $B \in \mathcal B_p(X, \tau)$. So $B = (O \setminus M) \cup N$, where $O$ is an open set in  $(X, \tau)$ and $M, N$ are meager sets. We can assume that $O \ne \emptyset.$
By Theorem~\ref{compatibility_theroem} $O = V \cup A$, where $V$ is a non-empty open set in 
$(X, \sigma)$ and $A$ is a nowhere dense set. Note that 
$B = ((V \cup A) \setminus M) \cup N = (V \setminus M) \cup ((A \setminus M) \cup N)$ and
$(A \setminus M) \cup N$ is a meager set. Hence, $B \in \mathcal B_p(X, \sigma)$, i. e. 
$\mathcal B_p(X, \tau) \subseteq \mathcal B_p(X, \sigma)$. The opposite inclusion can be proved 
similarly.
$\Box$

Besides the $\pi$- compatible topologies there are other answers to Question~\ref{question 1}.(b).
In fact,

\begin{example}\label{remark_baire} 
Let $\mathbb Q$ be the set of rational numbers, $\tau_E$ the Euclidean topology on $\mathbb Q$   
and $\tau_{\mbox{fc}}$ be the finite complement topology on $\mathbb Q$.
Since each one-point set of spaces $(\mathbb Q,\tau_E)$ and
$(\mathbb Q, \tau_{\mbox{fc}})$ is nowhere dense, by 
Proposition~\ref{countable_pro}  the families of their meager sets (resp. sets with the Baire property) are equal to $\mathcal P(\mathbb Q)$. 
 Note that $\tau_{\mbox{fc}}$ and $\tau_E$ are  not $\pi$-compatible.
\end{example}

\begin{question} Let $X$ be a set and $\tau, \sigma$ topologies on $X$.  What property can we add to the equality  $\mathcal{B}_p(X, \tau) = \mathcal{B}_p(X, \sigma)$
in order to get the $\pi$-compatibility of $\tau$ and $\sigma$?
\end{question}

Some weaker property than $\pi$-compatibility, the coincidence of the families of nowhere dense sets, does not garantee the equality between the families of the sets with the Baire property. In fact,

\begin{example} Let $\mathbb N$ be the set of positive integers, $\tau_{\mbox{dis}}$ be the discrete topology  and $\tau_p$ be the odd-even topology (see Example~\ref{topological sum}). Note that the families of nowhere dense sets (resp. meager sets) for the spaces $(\mathbb N, \tau_{\mbox{dis}})$ and $(\mathbb N, \tau_p)$ are equal to $\{\emptyset\}$.
But $\mathcal{B}_p(\mathbb N, \tau_p) =  \tau_p \ne \mathcal{B}_p(\mathbb N, \tau_{\mbox{dis}}) = \mathcal P(\mathbb N)$.

\end{example}

Recall that a topological space $X$ is called {\it a Baire space} if for every countable family $\{G_i\}_{i=1}^\infty$  of open dense subsets of $X$ the intersection $\cap_{i=1}^\infty G_i$ is dense in $X$.

It is evident that any open non-empty subset of a Baire space is also Baire.

\begin{theorem}\label{compatibility_Baire_space} Let $X$ be a set and $\tau, \sigma$ be $\pi$-compatible topologies on $X$.  Then the following is valid.
\begin{itemize}
\item[(a)] $(X, \tau)$ is a Baire space iff  $(X, \sigma)$ is  a Baire space.
\item[(b)] Let  $\mathcal G_\delta(X, \tau)$ (resp. $\mathcal G_\delta(X, \sigma)$ ) be the family of all $G_\delta$-sets of the space $(X, \tau)$ (resp. $(X, \sigma)$).

If $(X, \tau)$ is a Baire space (or  $(X, \sigma)$ is  a Baire space) then the family
$$\mathcal B = \{\emptyset \ne A \subseteq X : A \in \mathcal G_\delta(X, \tau) \cap \mathcal G_\delta(X, \sigma)\}$$
is a $\pi$-network for the space $(X, \tau)$ and for the space   $(X, \sigma)$.
\end{itemize}

\end{theorem}

Proof. (a) We show only the necessity.

Let $Y = \cup_{i=1}^\infty Y_i$, each $Y_i$ be closed in the space $(X, \sigma)$ and $\Int_{\sigma} Y_i = \emptyset$. By Lemma~\ref{compatibility_meager} we have that  $\Int_{\tau}\Cl_{\tau} Y_i = \emptyset$ for each $i$. Since the space  $(X, \tau)$ is Baire, we have $\Int_{\tau}(\cup_{i=1}^\infty \Cl_{\tau} Y_i) = \emptyset$. Assume that $O_\sigma = \Int_{\sigma}(\cup_{i=1}^\infty Y_i) \ne \emptyset$. Note that there is  $\emptyset \ne O_\tau \in \tau$ such that $O_\tau \subseteq O_\sigma \subseteq \cup_{i=1}^\infty Y_i.$ Hence, $O_\tau \subseteq \Int_{\tau}(\cup_{i=1}^\infty Y_i)$. Since $\emptyset \ne \Int_{\tau}(\cup_{i=1}^\infty Y_i) \subseteq \Int_{\tau}(\cup_{i=1}^\infty \Cl_{\tau} Y_i) = \emptyset$, we have a contradiction.

(b) Let us show that $\mathcal B$ is a $\pi$-network for the space $(X, \tau)$.

Consider a non-empty open set $O_1$ in $(X, \tau)$.  Apply Theorem~\ref{compatibility_theroem}. Choose a non-empty open set $V_1$ in
$(X, \sigma)$ such that $V_1 \subseteq O_1$ and $O_1 \setminus V_1$ is a nowhere dense set. Then choose a non-empty open $O_2$ in $(X, \tau)$ such that 
$O_2 \subseteq V_1$ and $V_1 \setminus O_2$ is a nowhere dense set, and so on. Since  
the space $(X, \tau)$ is a Baire space we have $\cap_{i=1}^\infty O_i \ne \emptyset$. Furthermore, by the construction
it follows that  $\cap_{i=1}^\infty O_i = \cap_{i=1}^\infty V_i \in \mathcal B$  and $\cap_{i=1}^\infty O_i$ is a dense subset of $O_1$. We have done.
$\Box$

\begin{proposition}\label{dense} Let $X$ be a set and $\tau, \sigma$ be $\pi$-compatible topologies on $X$. If a set $A \subset X$ is dense in $(X, \tau)$ then $A$ is dense in $(X, \sigma)$. In particular, $d(X, \tau) = d(X, \sigma).$ 
$\Box$
\end{proposition}

\vskip 0.3 cm

Different examples of $\pi$-compatible topologies one can get combining known example  and  the following simple statement.

\begin{proposition}\label{proposition_compatibility} (\cite[Proposition 3.2]{CN1}) Let $X_\alpha, \alpha \in \mathcal A,$ be  sets and  for each  $\alpha \in \mathcal A$  topologies $\tau_\alpha, \sigma_\alpha$ be $\pi$-compatible  on $X_\alpha$. Then
 the topological products of topologies $\prod \{\tau_\alpha, \alpha \in \mathcal A\}$ and 
$\prod \{\sigma_\alpha, \alpha \in \mathcal A\}$  on the set $\prod\{X_\alpha, \alpha \in \mathcal A\}$ are $\pi$-compatible.	$\Box$
\end{proposition}

Theorem~\ref{compatibility_Baire_space} and Propositions~\ref{dense}, ~\ref{proposition_compatibility} \ easily imply the following statement. 
\begin{corollary}\label{product_Baire} Let $\tau_i$ and  $\sigma_i$ be $\pi$-compatible topologies on a set $X_i $ for each $i \leq n$. Then 
\begin{itemize}
\item[(a)]  $d(\prod_{i=1}^n (X_i, \tau_i)) = d(\prod_{i=1}^n (X_i, \sigma_i)),$
\item[(b)] $\prod_{i=1}^n (X_i, \tau_i)$ is a Baire space iff  $\prod_{i=1}^n (X_i, \sigma_i)$ is  a Baire space. $\Box$
\end{itemize}

\end{corollary}

\section{Admissible extensions of topologies}

\begin{definition}\label{adm_ext} (\cite[Definition 3.1]{CN2})
A topology $\sigma$ on a set $X$ is said to be {\it an admissible extension of a topology $\tau$} on the same set $X$ if
$\tau \subseteq \sigma$; and 
$\tau$ is a $\pi$-base for $\sigma$, i.e. for each non-empty element $O$ of $\sigma$ there is a non-empty element $V$ of $\tau$ which is a subset of $O$. 
\end{definition}

Note that if $\tau$ and $\sigma$ are topologies on a set $X$ such that $\sigma$ is an admissible extension of $\tau$ then $\tau$ and $\sigma$ are $\pi$-compatible. 

\begin{theorem}\label{description_admissibility} Let $X$ be a set and $\tau, \sigma$ topologies on $X$ such that $\sigma$ is an admissible extension of $\tau$.
 If $O$ is a non-empty element of $\sigma$ then $O$ is a semi-open set of $(X, \tau)$, i. e. there is an element $V$ of $\tau$ such that $V \subseteq O \subseteq \Cl_{\tau} V$.
 \end{theorem}

Proof.  Put $V = \Int_{\tau} O$ and note that $V \ne \emptyset$. 
We will show that $\Cl_{\tau} V \supseteq O$. In fact,  assume that $W = O \setminus \Cl_{\tau} V \ne \emptyset.$ Since $\sigma$ is an admissible extension of $\tau$ then $W \in \sigma$ and there exists $\emptyset \ne U \in \tau$ such that $U \subseteq W \subseteq O$. It is easy to see that $U$ must be a subset of $V$. We have a contradiction which proves the statement. 
$\Box$

\begin{remark} Theorem~\ref{description_admissibility} does not hold for $\pi$-compatible topologies in general. Indeed, 
let $\tau_{S}$ be the Sorgenfrey  topology on the reals $\mathbb R$ generated by the family $\{[a,b): a, b \in \mathbb R\}$,
$\tau_{-S}$ be a topology on the reals $\mathbb R$ generated by the family $\{(a,b]: a, b \in \mathbb R\}$ and $\tau_E$ is the natural topology on $\mathbb R$. Note that both topologies $\tau_{S}$ and $\tau_{-S}$ are  admissible extensions of $\tau_E$. Hence, the topologies
$\tau_S$ and $\tau_{-S}$ are $\pi$-compatible. But for a non-empty element $(a,b] \in  \tau_{-S}$ we have the following:
$\Int_{\tau_S}(a,b] = (a,b)$ and $\Cl_{\tau_S}(a, b) = [a, b) \nsupseteq (a, b].$ $\Box$
\end{remark}

It is easy to see that the following proposition holds.

\begin{proposition}\label{axioms} Let $X$ be a set and $\tau, \sigma$ topologies on $X$ such that $\sigma$ is an admissible extension of $\tau$. 
If $(X, \tau)$ is a $T_0$-space (resp. a $T_1$-space, a $T_2$-space or a functionally Hausdorff space) then $(X, \sigma)$ is similar.  $\Box$ 

\end{proposition}

The inverse statement does not hold.

\begin{example}\label{double arrow} Let $X = X_0 \cup X_1 \subset \mathbb R^2$, where 
$X_i = \{\{x\} \times  \{i\} : x \in \mathbb R\}$. 

The topology $\sigma$ on $X$ is defined
by bases $\mathcal B_\sigma(p)$ at each point $p$ of $X$ as follows. 

If $p= \{x\} \times  \{1\}$ then $\mathcal B_\sigma(p) = \{[x,y) \times \{1\} \cup (x, y) \times \{0\}: y > x \}$.  

If $p= \{x\} \times  \{0\}$ then $\mathcal B_\sigma(p) = \{(y,x) \times \{1\} \cup (y, x] \times \{0\}: y < x \}$. 

The topology $\tau$ on $X$ is also defined
by bases $\mathcal B_\tau(p)$ at each point $p$ of $X$. Namely, 

If $p= \{x\} \times  \{1\}$ or  $p= \{x\} \times  \{0\}$ then 
$\mathcal B_\tau(p) = \{(y, z) \times \{1\} \cup (y, z) \times \{0\}: y < x < z \}$.

Let us note that the topology $\sigma$ is an admissible extension of the topology $\tau$ on the set $X$. Moreover, the space $(X, \sigma)$ is homeomorphic to a subspace of the Alexandroff double arrow space, and hence $(X, \sigma)$ is regular $T_1$, hereditarily Lindel\"of and hereditarily separable. However, $(X, \tau)$ is not $T_0$.

\end{example}

Proposition~\ref{axioms} cannot be extended to higher axioms.

\begin{example}\label{Smirnov} (\cite[Example 64]{SS}) Let us note that there exists an admissible extension $\sigma$ of the topology $\tau_E$ on the reals $\mathbb R$ such that $(\mathbb R, \sigma)$ is not regular.
As $\sigma$ we consider the Smirnov  topology on $\mathbb R$: $O \in \sigma$ if $O = U \setminus B$, where $U \in \tau_E$, $B \subseteq A$ and $A = \{\frac{1}{n} : n = 1,2, \dots\}$.
\end{example}

It is easy to see that if for some topologies $\sigma$ and $\nu$ on a set $X$ there exists a topology $\tau$ on $X$ such that
$\sigma$ and $\nu$ are admissible extensions of $\tau$ then $\nu$ and $\sigma$ are $\pi$-compatible, and $\tau \subseteq \sigma \cap \nu$.

\begin{remark} We will use the notations from Example~\ref{double arrow}. 

Let us consider the topology $\nu$ on $X$ defined
by bases $\mathcal B_\sigma(p)$ at each point $p$ of $X$ as follows. 

If $p= \{x\} \times  \{1\}$ then $\mathcal B_\sigma(p) = \{(y,x) \times \{0\} \cup (y, x] \times \{1\}: y < x \}$.  

If $p= \{x\} \times  \{0\}$ then  $\mathcal B_\sigma(p) = \{[x,y) \times \{0\} \cup (x, y) \times \{1\}: y > x \}$. 

Note that the topology $\nu$ is an admissible extension of the topology $\tau$, $\sigma \cap \nu = \tau$ and
the spaces $(X, \sigma)$ and $(X, \nu)$ are homeomorphic. 
\end{remark}

\begin{proposition}\label{topologies} Let $X$ be a set and $\tau, \sigma, \nu$ be topologies on $X$.
\begin{itemize}
\item[(a)] If $\sigma \subseteq \nu$ and $\sigma, \nu$ are admissible extensions of $\tau$ then $\nu$ is an admissble extension of $\sigma$.
\item[(b)] If $\tau \subseteq\sigma \subseteq \nu$ and $\nu$ ia an admissble extension of $\tau$ then 
$\sigma$ is an admissible extension of $\tau$. $\Box$
 \end{itemize}
 \end{proposition}
\begin{question} 
Let $X$ be a set, $\tau$ and $\sigma$ be $\pi$-compatible topologies on $X$. 
Are $\tau$ and $\sigma$  always admissible extensions of $\tau \cap \sigma$? 
\end{question}

\section{Some examples of admissible extensions}

\subsection{Topologies formed from a given topology and ideals of sets}

A way to obtain  finer topologies (but often non-regular) than a given topology $\tau$ on a set $X$ is to use ideals of subsets of $X$ (cf. \cite{JH}). 

Let  $\mathcal{I}$ be an ideal of subsets of  $X$ and let $A$ be a subset of $X$. For a point $x\in X$, let  $\mathcal{N}(x)=\{U\in \tau: x\in U\}$.

A local function of the set $A$ with the respect to the ideal $\mathcal{I}$  and the topology $\tau$, is the set $A^*(\mathcal{I},\tau)=\{x\in X: A\cap U\notin \mathcal{I} \text{ for every } U\in \mathcal{N}(x)\}$.

For every $A\subseteq X$ put $\Cl^{*}(A)=A\cup A^{*}(\mathcal{I}, \tau)$. 

It is easy to see that  $\Cl^{*}(\cdot)$ is a Kuratowski closure operator. So the family $\tau^{*}(\mathcal{I}) = \{U\subseteq X: \Cl^{*}(X\setminus U)=X\setminus U\}$ is a topology on the set $X$. 

Moreover, $\tau \subseteq \tau^{*}(\mathcal I)$, every element of the ideal $\mathcal I$ is closed in the space $(X, \tau^{*}(\mathcal I))$, and
$(\tau^{*}(\mathcal I))^{*}(\mathcal I) = \tau^{*}(\mathcal I)$.

It is easy to see that if $\mathcal I \subseteq \mathcal J$ then $\tau^{*}(\mathcal I) \subseteq \tau^{*}(\mathcal J)$.

\begin{theorem}\label{ideal_without_isolated_points} (\cite[Theorem 3.15]{CN1})
Let $(X, \tau)$ be a topological space without isolated points  and  $\mathcal{I}$  an ideal of subsets of 
$X$. Then the topology $\tau^*(\mathcal{I})$ is an admissible extension of $\tau$ if and only if $\mathcal{I} \subseteq \mathcal{I}_n(\tau),$ where $\mathcal{I}_n(\tau)$ is the family of nowhere dense sets of the space $(X, \tau)$. 

In particular,  $\tau^{*}(\mathcal I_n(\tau))$ is the finest admissible extension of $\tau$ formed from $\tau$ and  ideals of sets.
$\Box$
\end{theorem}

Applying Theorem~\ref{ideal_without_isolated_points} and Proposition~\ref{topologies}(a) we get
\begin{corollary} Let $(X, \tau)$ be a topological space without isolated points  and  $\mathcal{I}$,
$\mathcal{J}$  
 be ideals of subsets of 
$X$ such that $\mathcal{I} \subseteq \mathcal{J} \subseteq \mathcal{I}_n(\tau)$. Then 
$\tau^*(\mathcal{J})$ is an admissible extension of $\tau^*(\mathcal{I})$. $\Box$
\end{corollary}

\begin{remark} 
The condition on the space $(X, \tau)$ "to be without isolated points"  in Theorem~\ref{ideal_without_isolated_points} is essential. Indeed, let $(X, \tau_{\mbox{dis}})$ be a discrete space and $\mathcal{I} = \mathcal{P}(X)$. Note that  $ \tau_{\mbox{dis}}^*(\mathcal{I}) =  \tau_{\mbox{dis}}$, in particular,   $ \tau_{\mbox{dis}}^*(\mathcal{I})$ is an admissible extension of $ \tau_{\mbox{dis}}$, and $\mathcal{I}_n( \tau_{\mbox{dis}}) = \{\emptyset\}$.
\end{remark}

\begin{example}(cf. \cite{JH})
 Let $(X, \tau)$ be a topological space.
Then for each $A \subseteq X$ we have $A^*(\mathcal{I}_n(\tau),\tau)= \Cl_\tau (Int_\tau (\Cl_\tau (A)))$ and $\Cl^{*} (A) = A \cup \Cl_\tau (Int_\tau (\Cl_\tau (A)))$.  Moreover, $\tau^{*}(\mathcal I_n(\tau)) = \tau^\alpha$, where $\tau^\alpha$ is a topology from \cite{N}. 

\end{example}

Some facts about admissible extensions of the Euclidean topology $\tau_E$ on the real line  via ideals of sets see below. 

Since the real line $(\mathbb R, \tau_E)$ is a Baire space we get the following.

\begin{proposition} Let $\mathcal I$ be an ideal of sets on the real line $(\mathbb R, \tau_E)$ and $\mathcal I \subseteq \mathcal I_n(\tau_E)$. Then the space
 $(\mathbb R, \tau_E^{*}(\mathcal{I}))$ is connected.  $\Box$  
\end{proposition}

Let $\mathcal T_{\mbox{cd}}$ be the ideal of all closed discrete subsets of the real line $(\mathbb R, \tau_E)$. Recall (\cite[Proposition 4.2]{CN1}) that 
$\tau_E = \tau_E^*(\mathcal I)$ for any ideal $\mathcal I \subseteq \mathcal I_{\mbox{cd}}$.  

\begin{proposition} Let $\mathcal I$ be an ideal of sets on the real line $(\mathbb R, \tau_E)$ and $\mathcal I \subseteq \mathcal I_n(\tau_E)$. Then the space
 $(\mathbb R, \tau_E^{*}(\mathcal{I}))$ is non-regular iff there exists an element $I \in \mathcal I$ with a limit point in the topology $\tau_E$.  $\Box$  
\end{proposition}

\begin{remark}(cf. \cite{JH}) Let $A = \{\frac{1}{n} : n = 1,2, \dots\} \subset \mathbb R$ and $\mathcal{I}_A$ be the ideal of all subsets of $A$. Then the topology $\tau_E^{*}(\mathcal{I}_A)$ is the nonregular Smirnov topology on $\mathbb R$ from Example~\ref{Smirnov}.
\end{remark}

\subsection{Hattori spaces}
Other  extensions of the Euclidean topology $\tau_E$ on the real line 
was suggested in \cite{H}.

\begin{definition}\label{Hattori_definition}(\cite{H})
Let $A \subseteq \mathbb R$. Define a topology $\tau(A)$ on $\mathbb R$ as follows:
\begin{itemize}
	\item[(1)] for each $x \in A$, $\{(x-\epsilon, x+\epsilon): \epsilon > 0\}$ is a nbd open basis at $x$, 
	\item[(2)] for each $x \in \mathbb R \setminus A$, $\{[x, x+\epsilon): \epsilon > 0\}$ is a nbd open basis at $x$. 
\end{itemize}

\end{definition} 

Note that $\tau(\emptyset)$ (respectively, $\tau(\mathbb R)$) is the Sorgenfrey topology $\tau_S$ (respectively, the Euclidean topology $\tau_E$) on the reals, and
all topologies $\tau(A), A \subseteq \mathbb R,$ are addmissible extensions of $\tau_E$.

It is easy to see that for any $A, B \subseteq \mathbb R$ we have $A \supseteq B$ iff $\tau(A) \subseteq \tau(B)$.  

So by Proposition~\ref{topologies} we get

\begin{proposition}
$\tau(B)$ is an admissible extension of $\tau(A)$ iff $A \supseteq B$. $\Box$
\end{proposition}

The topological spaces $H(A) = (\mathbb R, \tau(A)), A \subseteq \mathbb R,$ are called {\it Hattori spaces}. 

Let us recall (\cite{CH}) that all $H(A)$-spaces are $T_1$ regular, hereditary Lindel\"of and hereditary separable. More on the properties of $H(A)$-spaces
and their variety see for example in \cite{K} and \cite{BS}. 

The following is evident.
 
\begin{proposition} 
Any $H(A)$-space with $A \ne \mathbb R$ is disconnected. $\Box$ 
\end{proposition}

\subsection{Almost topological groups}
{\it A paratopological group} $(G, \tau)$ consists of a group $G$ and a topology $\tau$ on $G$ that makes the group operation continuous. If in addition the inverse operation of $G$ is continuous then $(G, \tau)$ is {\it a topological group}.

{\it A semitopological group} is a group endowed with a topology that makes continuous left and right translations. Each paratopological group is a  semitopological one.

\begin{definition}(\cite[Definition 2.1]{F})\label{def_almost} {\it An almost topological group} is a paracompact group $(G, \tau)$ that satisfies the following conditions:
\begin{itemize}
\item[(a)] the group $G$ admits a Hausdorff topological group topology $\gamma$
weaker than $\tau$, and
\item[(b)] there exists a local base $\beta_e$ at the identity $e$ of the paratopological group $(G, \tau)$ such that the set $U \setminus \{e\}$ is open in $(G, \gamma)$ for every $U \in \beta_e$.
\end{itemize}
\end{definition}

One says also that $G$ is an almost topological group with structure $(\tau, \gamma, \beta_e)$.

It is easy to see that the topology $\tau$ is an admissible extension of $\gamma$.

\begin{proposition} (\cite[Proposition 2.5]{F}) Let $G$ be a non-discrete almost topological group with structure $(\tau, \gamma, \beta_e)$. Then $\gamma$ is the finest Hausdorff topological group topology on $G$ weaker than $\tau.$  
\end{proposition}

\begin{remark}(cf. \cite{CS}) For the paratopological group $(\mathbb R, \tau_S)$ 
the Euclidean topology  $\tau_E$ is the finest Hausdorff topological group topology on $\mathbb R$ weaker than $\tau_S.$  
\end{remark}

\begin{remark}
\begin{itemize}
\item[(a)]
The Hattori space $H(A)$, where $A$ is a proper subset of $\mathbb R$, is no semitopological group.
 
In fact, let $x \in \mathbb R \setminus A$ and $y \in A$. Note that $y + (x-y) = x$ and the translation of the additive group $\mathbb R$ by $x-y$ is not continuous at $y$.
\item[(b)]
However, the space $(\mathbb R, \tau_E^{*}(\mathcal I_n(\tau_E)))$ is a semitopological group with open shifts and continuous inversion but it is not a paratopological group. 

In fact, since the nowhere dense sets are invariant under translations and inversion as well as the open intervals we have that $(\mathbb R, \tau_E^{*}(\mathcal I_n(\tau_E)))$ is a semitopological group  with open shifts and continuous inversion.

Let us show that  $(\mathbb R, \tau_E^{*}(\mathcal I_n(\tau_E)))$ is no paratopological group.
Note that $0+0= 0$ and $W = (-1,1) \setminus \{\frac{1}{n} : n = 1,2, \dots\}$ is an open neighborhood of $0$
in $(\mathbb R, \tau_E^{*}(\mathcal I_n(\tau_E)))$. If  $(\mathbb R, \tau_E^{*}(\mathcal I_n(\tau_E)))$ is a paratopological group then there exist two open sets $V_1$ and $V_2$ containing $0$ 
such that $V_1 + V_2 \subseteq W$. Observe that we can find
a symmetric open neighborhood $V \subseteq V_1 \cap V_2$ of $0$  and a non-empty interval $(a, b)\subseteq V$. Since the difference $(a,b) - (a,b)$ contains a nondegenerated  interval $(-\epsilon, \epsilon)$
we have $(-\epsilon, \epsilon) \subseteq V - V \subseteq V_1 + V_2 \subseteq W$. But this is impossible.  
  \end{itemize}
\end{remark}

\begin{example}(\cite[Example 2.2]{F})\label{example} There exist non-regular almost topological groups. In fact, consider the additive group $\mathbb R^2$ with the identity
$e=(0,0)$. For every $r > 0$ let $B_r = \{e\} \cup \{(x, y) \in \mathbb R^2 : x^2+y^2 < r^2, x > 0\}$. It is easy to see that the family $\mathcal B = \{B_r : r> 0\}$ is a local base at $e$ of a Hausdorff paratopological group. 
The group $(\mathbb R^2, \tau)$, where $\tau$ is a topology generated by the family $\mathcal B$, is a Hausdorff non-regular almost topological group. 
\end{example}

\subsection{Hattori topologies on almost topological groups}

We will follow (\cite{CS}) and apply notations  from Definition~\ref{def_almost}. 

Let $(G, \tau)$ be an almost topological group and $\beta_e$ a local base at the identity $e$ of $(G, \tau)$. Then $\{UU^{-1} : U \in \beta_e\}$ is a local base at $e$ in the Hausdorff topological group $(G, \gamma)$. Given $x \in G$, let $\beta_x = \{ Ux : U \in \beta_e \}$ and $\beta_e' = \{UU^{-1}x : U \in \beta_e\}$.

Given $A \subseteq G$, consider the collection $\{\beta(x) : x \in G\}$, where
$$\beta(x) =  \begin{array}{l} \beta_x'. \ \mbox{if}  \ x \in A, \\ 
\beta_x. \ \mbox{if}  \ x \notin A. \end{array}$$

\begin{theorem}(\cite[Theorem 3.1]{CS}) Let $(G, \tau)$ be an almost topological group and $A \subseteq G$. Then $\{\beta(x) : x \in G\}$ is a neighborhood system for a topology $\tau(A)$ on $G$. $\Box$
\end{theorem}
The topological space $(G, \tau(A))$ is denoted by $H(A,G)$ and it is called {\it the Hattori space} of $G$ associated to $A$. Let us note that $H(\emptyset,G) = (G, \tau)$ and  $H(G,G) = (G, \gamma)$.

\begin{proposition}\label{implication}(\cite[Proposition 3.6]{CS}) Let 
 $(G, \tau)$ be an almost topological group such that $(G, \tau)$ is not a topological group.
 If $A, B \subseteq G$ then $\tau(A) \subseteq \tau(B)$ iff $B \subseteq A$.
 $\Box$
\end{proposition}

Proposition~\ref{implication} and Proposition~\ref{topologies} imply
\begin{proposition} 
Let 
 $(G, \tau)$ be an almost topological group such that $G$ is not a topological group, and $A, B\subseteq G$.
 Then $\tau(A)$ is an admissible extension of $\tau(B)$  iff $A \subseteq B$.
 $\Box$
\end{proposition}

\begin{problem}  Let $(G, \tau)$ be the almost topological group and $\gamma$ be the finest Hausdorff topological group topology on $G$ weaker than $\tau.$ (For example, let $(G, \tau)$ be the almost topological group from Example~\ref{example}.) Describe the topological diversity of Hattori spaces   
$H(A,G), A \subseteq G$.
\end{problem}

\noindent(V.A. Chatyrko)\\
Department of Mathematics, Linkoping University, 581 83 Linkoping, Sweden.\\
vitalij.tjatyrko@liu.se

\end{document}